\documentclass[a4paper, 12pt]{article}

\usepackage[koi8-r,cp1251]{inputenc}
\usepackage[english,russian]{babel}
\usepackage{amsfonts,amssymb,amsmath, hyperref}
\usepackage[final]{epsfig}
\usepackage{graphicx}

\begin{document}

\begin{center}
 A. R. Mirotin\\
\textbf{FREDHOLM AND SPECTRAL PROPERTIES OF  TOEPLITZ OPERATORS ON $H^p$ SPACES OVER ORDERED GROUPS}\\
 amirotin@yandex.ru
\end{center}

Toeplitz operators on  spaces $H^p(G)\ (1< p<\infty)$ associated with
compact  connected Abelian group $G$ with ordered dual are considered and the generalization
of the classical Gohberg-Krein theorem on the Fredholm index of such  operators
with continuous symbols is proved.  Applications to spectral theory
of Toeplitz operators are given and  examples of evident computation
of index  have been considered.

%\emph{Key words and phrases: Toeplitz operator, Fredholm operator, Fredholm index,
%essential spectrum, odered abelian group}.

\

%УДК 517.983.23+517.984.5

\begin{center}
\bf А. Р. Миротин  \\
ФРЕДГОЛЬМОВЫ И СПЕКТРАЛЬНЫЕ СВОЙСТВА ТЁПЛИЦЕВЫХ ОПЕРАТОРОВ В ПРОСТРАНСТВАХ $H^p$
НАД УПОРЯДОЧЕННЫМИ ГРУППАМИ
%\footnote[1]{Работа выполнена при финансовой поддержке
%Государственной программы фундаментальных исследований Республики
%Беларусь, договор \No~20061473.}
\end{center}

{\small Рассматриваются тёплицевы
операторы в пространствах $H^p(G)\ (1< p<\infty)$, ассоциированных с
компактной связной абелевой группой $G$,  группа характеров которой   упорядочена, и в случае
линейного порядка доказывается теорема об индексе Фредгольма для
таких операторов с непрерывным символом, обобщающая классическую
теорему Гохберга-Крейна.  Указываются приложения
 полученных результатов к спектральной теории тёплицевых операторов и
 рассматриваются примеры явного вычисления индекса.
 }

%\emph{Ключевые слова: тёплицев оператор, фредгольмов оператор, индекс Фредгольма,
%существенный спектр, упорядоченная абелева группа.}

\

\vspace{8mm}
\centerline{\bf \S 1.  Введение}
\vspace{5mm}

Классическая теория тёплицевых операторов, возникшая первоначально
как теория операторов в гильбертовом пространстве $H^2$ на
группе $\mathbb{T}$ вращений окружности, была в дальнейшем перенесена
на банаховы пространства $H^p\ (1<p<\infty)$ на  $\mathbb{T}$ (см., например,  \cite{BS}).
Данная работа посвящена некоторым вопросам фредгольмовой и спектральной
теории тёплицевых операторов в пространствах $H^p\ (1<p<\infty)$
на компактных абелевых группах.

 Всюду ниже   $G$ есть нетривиальная связная
компактная абелева группа с  нормированной мерой Хаара $m$ и
(вообще говоря, частично) упорядоченной группой характеров $X$,
$T$ -- положительный конус в  $X$, порождающий $X$. Другими
словами, в группе  $X$ выделена   подполугруппа $T$, содержащая
единичный характер ${\bf 1}$
и такая, что $T\cap T^{-1}=\{{\bf 1}\}$ и $X=T^{-1}T$. При этом
полугруппа $T$ индуцирует в $X$  частичный порядок, согласованный
со структурой группы,  по правилу $\xi\leq\chi:=\chi\xi^{-1}\in
T$. Очевидно, что $T$ -- конфинальная часть $X$ и направленное
множество относительно рассматриваемого порядка (если  $\lambda_1,\lambda_2\in T$,
то $\lambda_1\leq \lambda_1\lambda_2,\lambda_2\leq \lambda_1\lambda_2$). Если
дополнительно $T\cup T^{-1}=X$, то этот порядок линеен (пишем $X$
л. у.). Хорошо известно, что (дискретная) абелева группа $X$ может быть линейно
упорядочена тогда и только тогда, когда она не имеет кручения
(см., например, \cite{Rud}), что, в свою очередь, равносильно
тому, что её группа характеров $G$ связна \cite{Pont} (при этом линейный порядок
в $X$, вообще говоря, не единственен).

В описанной выше ситуации
 тёплицевы операторы
$T_\varphi$ с символом $\varphi\in L^\infty(G)$ в обобщённых
пространствах Харди $H^2(G)$  были определены Дж. Мёрфи в \cite{Mur87} и
интенсивно изучались (см., например, \cite{Mur87} -- \cite{XuChen}). В частности, было показано,
что соответствующие обобщённые тёплицевы алгебры обладают рядом интересных свойств,
если порядок линеен. Более общее определение тёплицева оператора над группами дано
в \cite{IJPAM}, \cite{Corr}.
Рассматривались также обобщения и дискретных тёплицевых операторов, в том числе и
с матричным символом,
см. статью \cite{EMRS} и библиографию там (автор признателен рецензенту, обратившему
его внимание на эту работу).

Ниже мы рассматриваем тёплицевы
операторы в пространствах $H^p(G)\ (1< p<\infty)$ и в случае
линейного порядка доказываем теорему об индексе Фредгольма для
таких операторов с непрерывным символом, обобщающую классическую
теорему Гохберга-Крейна. Попутно установлены обобщения
теорем Брауна-Халмоша и Хартмана-Винтнера.  Указаны  приложения
 полученных результатов к спектральной теории тёплицевых операторов и рассмотрены примеры
 явного вычисления индекса.
 При этом, хотя ряд результатов для общих групп $G$ оказался
 аналогичен классическому случаю $G=\mathbb{T}$, обнаружились и существенные отличия.
 Например, если $G\not=\mathbb{T}$, то не каждый тёплицев оператор с
 неаннулирующимся непрерывным символом фредгольмов,
 существуют группы, для которых
 фредгольмовость соответствующих им тёплицевых операторов равносильна
 обратимости, существенный спектр тёплицева оператора,
 вообще говоря, не совпадает с множеством значений его символа
и т. д.  (см. теорему 4 и следствия 1, 2, 3, 5, 6, 7 и 11, а также примечание рецензента
к примеру 2). Таким образом,
 в определенном смысле группа  $\mathbb{T}$  в этом круге
 вопросов оказалась уникальной.

 Отметим, что в
 работах \cite{MurIrish91} --- \cite{Mur93} изучались
 индексы и обобщённые индексы тёплицевых операторов в
пространствах $H^2(G)$ (случай
 архимедова порядка рассматривался ранее в \cite{CDSS}).
 В частности, в \cite{Mur93} получено обобщение
теоремы об индексе при условии, что группа $X$ линейно упорядочена и содержит наименьший
положительный элемент. При этом использовались методы теории
$C^\ast$-алгебр, не переносящиеся на случай банаховых пространств
$H^p(G)$.

\vspace{5mm}
\centerline{\bf \S 2.  Тёплицевы операторы над упорядоченными группами}
\vspace{5mm}

Ниже через $Pol(G)$ ($Pol_T(G)$) будет обозначаться пространство
тригонометрических полиномов (соответственно, тригонометрических
полиномов аналитического типа) на группе $G$, т. е. линейная
оболочка множества $X$  (соответственно, $T$) в $L^p(G)$. При
$1\leq p\leq\infty$ через $H^p(G)$ обозначим подпространство таких
функций $f \in L^p(G)$, преобразование Фурье $\widehat f$ которых
сосредоточено на $T$, с нормой, индуцированной из $L^p(G)$
\cite{HL}.  Легко видеть, что $Pol_T(G)$  содержится в  $H^p(G)$ (кроме того,
$A_T(G)\subset H^\infty(G)$, где алгебра $A_T(G)$ есть равномерное
замыкание  множества $Pol_T(G)$, см. лемму 2).

%\begin{lemma}
\textbf{Лемма 1.} \textit{При $1\leq p<\infty$ пространство $H^p(G)$ инвариантно относительно
преобразования  $g\mapsto g^\ast$, где $g^\ast(y)=\overline {g(y^{-1})}$, и совпадает с
 замыканием  в $L^p(G)$ множества $Pol_T(G)$.}
%\end{lemma}

Доказательство. Инвариантность $H^p(G)$ относительно преобразования $g\mapsto g^\ast$
проверяется непосредственно с учётом инвариантности меры Хаара группы $G$ относительно отражения.
 Обозначим замыкание в $L^p(G)$ множества $Pol_T(G)$ через $H^p_T(G)$.
Включение $H^p_T(G)\subseteq H^p(G)$  следует из того, что
$H^p(G)$ замкнуто в  $L^p(G)$ и содержит   $Pol_T(G)$.

Пусть  $F$ есть такой линейный непрерывный функционал на $L^p(G)$,
что сужение $F|H^p(G)=0$. Тогда найдётся  функция $f\in L^q(G)$  (здесь и ниже $p^{-1}+q^{-1}=1$),
для которой
$$
F(h)=\int\limits_G h(x)\overline{f(x)}dm(x)=0
$$
\noindent при всех $h\in  H^p(G)$.
В частности, преобразование Фурье $\widehat f(\chi)=0$  при всех $\chi\in T$.  Если
теперь $g\in H^p(G)$, то   $\widehat f\widehat g=0$, а тогда и свёртка $f\ast g=0$, т. е.

$$
\int\limits_G f(xy^{-1})g(y)dm(y)=0.
$$
\noindent при почти всех $x\in G$.  Отсюда следует, что $F(g^\ast)=0$.
Следовательно, по теореме Хана-Банаха   $g^\ast\in H^p_T(G)$, а потому и
$g\in H^p_T(G)$, поскольку последнее пространство также инвариантно
относительно преобразования $g\mapsto g^\ast$.

%\begin{lemma}
\textbf{Лемма 2. }{\rm 1)} \textit{Для любого $p\in [1,\infty]$ имеет место включение $H^\infty(G)\cdot H^p(G)\subseteq H^p(G)$.
Более того, для каждого фиксированного $p\in [1,\infty)$}
$$
H^\infty(G)=\{h\in L^\infty(G)| hH^p(G)\subseteq H^p(G)\}.\eqno(2.1)
$$

{\rm 2)}\textit{ Пространство $H^\infty(G)$ есть банахова алгебра относительно поточечных операций и нормы,
индуцированной из  $L^\infty(G)$}.
%\end{lemma}

Доказательство. 1) Пусть $h\in H^\infty(G),\ f\in H^p(G), p\in [1,\infty)$. Тогда $hf\in L^p(G)$,
и при $\chi\in X\setminus T,\ \xi\in T$ получаем $\widehat{h\xi}(\chi)=\widehat{h}(\xi^{-1}\chi)=0$,
так как $\xi^{-1}\chi\in X\setminus T$. Следовательно, $\widehat{hg}(\chi)=0$ для
любого аналитического полинома  $g\in Pol_T(G)$. Выберем теперь последовательность $g_n\in Pol_T(G)$,
сходящуюся к $f$ в  $L^p(G)$ (лемма 1). Если $\chi\in X\setminus T$, то по доказанному выше
с учётом неравенства Гёльдера ($h\overline{\chi}\in L^q(G)$)
$$
0=\widehat{hg_n}(\chi)=\int\limits_G h\overline{\chi}g_ndm\to\widehat{hf}(\chi)\ (n\to\infty),
$$
а потому $hf\in H^p(G)$.

Если же $h\in L^\infty(G)$ и $hH^p(G)\subseteq H^p(G)$, то $h\cdot 1\in H^p(G)$, что доказывает
(2.1). В свою очередь,  (2.1) влечёт включение $H^\infty(G)\cdot H^\infty(G)\subseteq H^\infty(G)$.

 2) Это следует из 1) и замкнутости  $H^\infty(G)$ в $L^\infty(G)$.

%\begin{definition}
\textbf{Определение 1. } {\it Проектор Рисса} $P_T: Pol(G)\to Pol_T(G)$ определяется
равенством
$$
P_T(\sum\limits_{\chi\in M}c_\chi \chi)=\sum\limits_{\chi\in M\cap
T}c_\chi \chi.
$$
%\end{definition}

Известно  \cite{B}, \cite{Hel} (см. также \cite[Глава 8]{Rud}), что $P_T$
$L^p$-ограничен при $1<p<\infty$, а потому продолжается до
ограниченного проектора $P_T:L^p(G)\to H^p(G)$ (теорема Бохнера-Хелсона).

%\begin{definition}
\textbf{Определение 2. } Пусть $1< p<\infty$. Тёплицев оператор $T_\varphi$ в $H^p(G)$ с
символом $\varphi\in L^\infty(G)$  определяется  следующим образом:
$$
T_\varphi f=P_T(\varphi f),\ f\in H^p(G).
$$
%\end{definition}

Далее через $M_\varphi$ будет обозначаться оператор умножения на измеримую функцию
$\varphi$, действующий в пространстве $L^p(G)$ (или его подпространствах).

 Нижеследующее обобщение теоремы Брауна-Халмоша  (относительно классического
 случая см., например, \cite{BS}) справедливо для
произвольной компактной абелевой группы $G$  с линейно
упорядоченной группой характеров. Для его доказательства нам понадобится

%\begin{lemma}
\textbf{Лемма 3. } \textit{Если для измеримой функции $\varphi$ на $G$ оператор умножения
$M_\varphi$ отображает $Pol(G)$ в $L^p(G)$ и $L^p$-ограничен ($1<p<\infty$),
то $\varphi\in L^\infty (G)$, и  $\|M_\varphi\|=\|\varphi\|_\infty$.}
%\end{lemma}

Доказательство. Интерес представляет случай $M_\varphi\ne O$. Допустим, что
$\varphi\notin L^\infty (G)$.
Тогда найдётся такое
компактное $K\subset G$, что $m(K)>0$ и $|\varphi(x)|>2\|M_\varphi\|$ при $x\in K$.
Для любого $\varepsilon>0$ выберем  открытое $U\supset K$, удовлетворяющее
условию  $m(U)<m(K)+\varepsilon$,
и пусть непрерывная функция $f:G\to [0;1]$ такова, что $f|K=1,\ f|(G\setminus U)=0$
(лемма Урысона).
По теореме Вейерштрасса-Стоуна существует  полином $q\in Pol(G)$, такой, что
$\|f-q\|_\infty<\varepsilon$. Отсюда следует, что $|q(x)|\geq 1-\varepsilon$ при $x\in K$, а потому
$$
\|M_\varphi q\|_p\geq \left(\int\limits_K|\varphi|^p|q|^p dm\right)^{1/p}\geq
(1-\varepsilon)\left(\int\limits_K|\varphi|^p dm\right)^{1/p}>
$$
$$
>(1-\varepsilon)2\|M_\varphi\|(m(K))^{1/p}.
$$

С другой стороны, $\|f-q\|_p<\varepsilon$, и значит $\|q\|_p<\|f\|_p+\varepsilon$.
А поскольку к тому же $\|f\|_p<m(U)^{1/p}$, то
$$
\|M_\varphi q\|_p\leq \|M_\varphi\|(\|f\|_p+\varepsilon)\leq \|M_\varphi\|(m(U)^{1/p}+\varepsilon)<
$$
$$
<\|M_\varphi\|((m(K)+\varepsilon)^{1/p}+\varepsilon).
$$

Следовательно, для любого $\varepsilon>0$ справедливо неравенство
$$
(1-\varepsilon)2\|M_\varphi\|m(K)^{1/p}<\|M_\varphi\|((m(K)+\varepsilon)^{1/p}+\varepsilon),
$$
что невозможно.
Последнее утверждение леммы хорошо известно.

Далее мы полагаем для  $f\in   L^p(G), g\in L^q(G)
(p\in(1,\infty), p^{-1}+q^{-1}=1)$
$$
 \langle f,g \rangle:= \int\limits_G f\overline g dm.
$$

 %\begin{theorem}
\textbf{Теорема 1.}  \textit{Пусть $X$  л. у.,  $1< p<\infty$.
 Если ограниченный оператор $A:Pol_T(G)\to
H^p(G)$ таков, что для некоторой функции $a:X\to \mathbb{C}$ при любых $\chi_1, \chi_2\in T$
справедливо равенство  $\langle A\chi_1, \chi_2\rangle=a(\chi_1{ \chi_2}^{-1})$, то существует такая функция
$\varphi\in L^\infty(G)$, что $A=T_\varphi$, причём её преобразование Фурье $\widehat\varphi=a$. Более того,
 $ \|\varphi\|_\infty\leq \|T_\varphi\|\leq
c_p\|\varphi\|_\infty$, где $c_p=\|P_T\|$. В частности, $\|T_\varphi\|
=\|\varphi\|_\infty$, если $p=2$.}
%\end{theorem}

Доказательство. Рассмотрим направленность $b_\chi=\overline{\chi}A\chi\in L^p(G)\ (\chi\in T)$.
Так как она ограничена ($\|b_\chi\|_p\leq \|A\|$ для всех $\chi\in T$), то, переходя, если нужно, к
поднаправленности, можно считать, что $b_\chi$ слабо сходится в $L^p(G)$ к элементу $\varphi\in L^p(G)$.
В частности, $\lim_{\chi\in T}\langle b_\chi,\xi\rangle=\langle \varphi,\xi\rangle$ при всех $\xi\in X$.
Заметим, что $\langle b_\chi,\xi\rangle=\langle A\chi,\chi\xi\rangle=a(\xi^{-1})$,
если $\chi\xi\in T$. Поскольку $T$ конфинально в  $X$,  в пределе получаем, что
$\langle \varphi,\xi\rangle=a(\xi^{-1})$ при всех $\xi\in X$, т. е. $\widehat\varphi=a$.

Оператор $M_\varphi$,  очевидно, отображает $Pol(G)$ в $L^p(G)$;
покажем, что он $L^p$-ограничен. В силу последнего равенства имеем при всех $f,g\in Pol(G)$
$$
\langle M_\varphi f,g\rangle=\langle M_{\bar{\chi}}AM_\chi f,g \rangle,
$$
так как обе части совпадают при $f,g\in X$. Поэтому с учётом неравенства Гёльдера
$$
|\langle M_\varphi f,g\rangle|\leq\| M_{\bar{\chi}}AM_\chi f\|_p\|g\|_q\leq \|A\|\|f\|_p\|g\|_q,
$$
и следовательно
$$
\|M_\varphi f\|_p=\sup\{|\langle M_\varphi f,g\rangle|| g\in Pol(G), \|g\|_q\leq 1\}\leq \|A\|\|f\|_p.
$$
Таким образом, оператор $M_\varphi$ $L^p$-ограничен, а потому
$\varphi\in L^\infty (G)$ по лемме 3. Одновременно мы показали, что  $\|M_\varphi\|\leq \|A\|$.

Далее, так как при $\chi_1,\chi_2\in T$
$$
\langle T_\varphi \chi_1,\chi_2\rangle=\widehat\varphi(\chi_1\chi_2^{-1})=a(\chi_1\chi_2^{-1})
=\langle A\chi_1,\chi_2\rangle,
$$
то $A=T_\varphi$.

Наконец,
$$
\|\varphi\|_\infty=\|M_\varphi\|\leq\|A\|=\|T_\varphi\|=\|P_T M_\varphi\|\leq c_p\|\varphi\|_\infty,
$$
и доказательство теоремы полностью завершено.

\vspace{5mm}
\centerline{\bf \S 3. Фредгольмовость и индекс тёплицевых операторов}
\centerline {\bf с непрерывным символом}
\vspace{5mm}

Для формулировки и доказательства теоремы об индексе требуется определенная
подготовка.

Прежде всего, мы определим индекс вращения для функций из некоторой
подгруппы группы $C(G)^{-1}$ обратимых элементов алгебры $C(G)$.
Начнём с определения индекса вращения характера группы $G$ (ниже
$\#F$ обозначает число элементов конечного множества $F$).

 %\begin{definition}
\textbf{Определение 3.}  В каждом из следующих случаев определим индекс  вращения характера $\chi\in X$
 следующим образом:

 {\rm 1)} ${\rm ind}\chi = \#(T\setminus \chi T)$, если $\chi\in T$ и множество
 $T\setminus \chi T$ конечно;

{\rm 2)}  ${\rm ind}\chi ={\rm ind}\chi_1-{\rm ind}\chi_2$, если
   $\chi=\chi_1\chi_2^{-1}$, где $\chi_j\in T$, причём оба множества
 $T\setminus \chi T_j$  конечны  $(j=1,2)$.

  В остальных случаях считаем, что характер не имеет индекса.

%\end{definition}

Множество характеров, имеющих индекс, обозначим $X^i$.

Далее нам понадобится результат Г.~Бора и E.~ван~Кампена, согласно
которому любая функция $\varphi\in C(G)^{-1}$ представима в виде
$\chi e^g$, где $g\in C(G),\ \chi\in X$ \cite{vK}. Более того, как
следует из одного результата  Е.~А.~Горина (см.  замечание после доказательства теоремы 2 в \cite{Gor}),
характер $\chi$  в этом разложении определяется по $\varphi$
однозначно (позднее этот факт был переоткрыт в \cite{MurIrish91}).

%\begin{definition}
\textbf{Определение 4.} Рассмотрим функцию $\varphi\in C(G)^{-1}$ с разложением
Бора-ван Кампена  $\varphi=\chi e^g\ (g\in C(G),\ \chi\in X)$.
 Если  $\chi\in
X^i$, то положим
$$
{\rm ind}\varphi={\rm ind}\chi.
$$
В противном случае будем считать, что функция $\varphi$ не имеет индекса.
%\end{definition}

Множество функций из $C(G)^{-1}$, имеющих индекс, обозначим
$\Phi(G)$. Таким образом, $\Phi(G)=X^i\cdot \exp(C(G))$, причём по
причине отмеченной выше единственности характера в разложении
Бора-ван Кампена, $X\cap \exp(C(G))=\{{\bf 1}\}$.

В следующей теореме перечислены основные свойства множеств $X^i$, $\Phi(G)$
и отображения ${\rm ind}\ (H^\bot$ обозначает аннулятор подмножества $H$  группы
 $G$, а  знак $\sqcup$ --- дизъюнктное объединение множеств;
 напомним, что подгруппа $\Xi$ упорядоченной группы $X$ называется \textit{выпуклой}, если отрезок
 $[\xi_1,\xi_2]$ содержится в  $\Xi$, лишь только $\xi_1,\xi_2\in\Xi$).

%\begin{theorem}
\textbf{Теорема 2.} {\rm 1)} \textit{Множество $X^i$ есть выпуклая подгруппа группы }$X$;

{\rm 2)} \textit{отображение ${\rm ind}:X^i\to \Bbb{Z}$ определено
корректно и есть сохраняющий порядок гомоморфизм на группу  вида $a\Bbb{Z}
\ (a\in \Bbb{Z}_+)$, который, если $X$ л. у., а  $X^i$ нетривиальна,
является порядковым изоморфизмом на } $\Bbb{Z}$;

{\rm 3)} \textit{множество $\Phi(G)$ есть открытая подгруппа группы
$C(G)^{-1}$ с открытыми компонентами связности  $\chi\exp(C(G))$,
где $\chi$ пробегает $X^i$;
}
{\rm 4)} \textit{отображение  ${\rm ind}:\Phi(G)\to \Bbb{Z}$  есть
гомоморфизм групп и принимает постоянное значение ${\rm ind}\chi$
на компоненте связности  $\chi\exp(C(G))$ ($\chi\in X^i$);
}
{\rm 5)} \textit{если $X$ л. у., то для любого непрерывного гомоморфизма
$\gamma:\mathbb{T}\to G$, удовлетворяющего условию $X^i\nsubseteq\gamma(\mathbb{T})^\bot$,
найдётся такое целое $k_\gamma\ne 0$, что при всех $\varphi\in \Phi(G)$ справедливо равенство}
$$
{\rm ind}\varphi=k_\gamma^{-1}{\rm wn}(\varphi\circ\gamma),
$$
\textit{где ${\rm wn}$ обозначает классический индекс вращения функций из $C(\mathbb{T})^{-1}$.
}%\end{theorem}

Доказательство. 1), 2) Сначала заметим, что
$$
  T\cap X^i=\{\chi\in T|  \#(T\setminus \chi T)<\infty\}.
$$
Действительно, если $\chi\in T\cap X^i$, то
$\chi=\chi_1^{-1}\chi_2$, причём
 $\#(T\setminus \chi_j T)<\infty\  (j=1,2)$. Тогда  $\chi_2=\chi_1\chi$
 (т. е. ${\bf 1}\leq\chi\leq\chi_2$), а потому  $T\setminus \chi T\subseteq T\setminus \chi_2
 T$. Отсюда следует, что   $\#(T\setminus \chi T)<\infty$, и ${\rm
 ind}\chi\leq {\rm ind}\chi_2$   (заодно мы показали, что $\chi_2\in T\cap X^i$ влечёт
 $[{\bf 1},\chi_2]\subset X^i$, и что отображение ${\rm ind}:T\to \Bbb{Z}_+$ сохраняет порядок).
Значит, $T\cap X^i\subseteq\{\chi\in T|  \#(T\setminus \chi T)<\infty\}.$ Обратное включение очевидно.

Установим корректность определения ${\rm ind}\chi$. При $\chi_1,
\chi_2 \in T$           справедливо равенство
$$
T\setminus \chi_1\chi_2 T=(T\setminus \chi_1
T)\sqcup\chi_1(T\setminus \chi_2 T).
$$
Поэтому, если $\chi_1, \chi_2$ имеют индексы, то $\chi_1 \chi_2$
тоже его имеет (т. е. $T\cap X^i$ есть полугруппа) и   ${\rm
ind}(\chi_ 1\chi_2)={\rm ind}\chi_1+{\rm ind}\chi_2$. Из
последнего равенства следуют  корректность определения ${\rm
ind}\chi$ (включая  совпадение  двух определений ${\rm
ind}\chi$ в случае $\chi\in T$) и гомоморфность индекса.

Далее, $X^i=(T\cap X^i)^{-1}(T\cap X^i)$    есть подгруппа группы $X$, порождённая
полугруппой  $T\cap X^i$.
Докажем её выпуклость.    Пусть $\xi_1\leq\xi\leq\xi_2$, где    $\xi_1=\chi_1^{-1}\chi_2,\
 \xi_2=\chi_3^{-1}\chi_4,\ (\chi_j\in T\cap X^i, j=1,\ldots,4)$.
 Тогда $\chi_2\chi_3\leq\chi_1\chi_3\xi\leq\chi_1\chi_4$, причём из
 первого неравенства следует, что $\chi_5:=\chi_1\chi_3\xi\in T$,
 а из второго, -- что $\chi_5\in X^i$ (см. начало доказательства
 теоремы). Поэтому $\xi=(\chi_2\chi_3)^{-1}\chi_5\in X^i$.

 Покажем, что отображение  ${\rm ind}$ сохраняет порядок. Пусть,
 как выше,  $\xi_1\leq\xi\leq\xi_2$, где    $\xi_1=\chi_1^{-1}\chi_2,\
 \xi_2=\chi_3^{-1}\chi_4,\ (\chi_j\in T\cap X^i, j=1,\ldots,4)$.
Тогда $\chi:=\xi_2\xi_1^{-1}\in T\cap X^i$. Поэтому $\chi_2\chi_3\chi=\chi_1\chi_4$,
 откуда следует, что  ${\rm ind}\chi_2+{\rm ind}\chi_3\leq {\rm ind}\chi_1+{\rm ind}\chi_4$,
  т. е.
 ${\rm ind}\xi_1\leq {\rm ind}\xi_2$. Кроме того, ${\rm ind}(X^i)$,
 будучи подгруппой группы $\Bbb{Z}$, имеет вид $a\Bbb{Z}\ (a\in\Bbb{Z}_+)$.

 Предположим теперь, что   $X$ л. у., а  $X^i$ нетривиальна, и докажем инъективность
 гомоморфизма  ${\rm ind}:X^i\to \Bbb{Z}$, т. е. равенство ${\rm
 Ker}({\rm ind})=\{{\bf 1}\}$.     Пусть $\xi=\chi_1^{-1}\chi_2,\ \chi_j\in T\cap X^i, j=1,2$
 и  ${\rm ind}\xi=0$, т. е.  ${\rm ind}\chi_1={\rm ind}\chi_2$.
 Если, для определённости, $\chi_1\leq\chi_2$, то   $T\setminus \chi_1 T
 \subseteq T\setminus \chi_2 T$, а потому здесь имеет место
 равенство (это конечные множества с одинаковым числом элементов), откуда
 $\chi_1 T=\chi_2 T$.   Последнее равенство влечёт
 $\chi_1\leq\chi_2$ и   $\chi_2\leq\chi_1$, т. е.  $\chi_1=\chi_2$.
 Таким образом, группа $X^i$ порядково изоморфна $a\Bbb{Z}$, а потому
  ${\rm ind}(X^i)={\rm ind}(a\Bbb{Z})=\Bbb{Z}$, поскольку  $a\not=0$.

 3)     Равенство $\Phi(G)=X^i\cdot \exp(C(G))$ показывает,
 что   $\Phi(G)$ есть  подгруппа группы      $C(G)^{-1}$ (причём из равенства
 $X^i\cap \exp(C(G))=\{{\bf 1}\}$
 следует, что она изоморфна прямому произведению $X^i\times \exp(C(G))$).
 Далее, хорошо известно, что   множество  $\exp(C(G))$  открыто и является
 компонентой единицы  в      $C(G)^{-1}$.  Следовательно,
 подмножества  $\chi\exp(C(G))\subset C(G)^{-1}$ открыты и связны,
 и осталось заметить, что при различных $\chi$    они попарно не
 пересекаются   по причине  единственности характера в разложении
Бора-ван Кампена.

4)  Это следует из 2) и того факта, что   ${\rm ind}(\exp(C(G)))=\{ 1\}$.

5) Так как ${\rm ind}(e^g)={\rm wn}(e^{g\circ\gamma})=0$ при $g\in C(G)$, интерес представляет случай
$\varphi\in X^i$. С учётом утверждения 2) можно, исключая тривиальный случай,
 считать также, что  группа $X^i$ изоморфна $\mathbb{Z}$. Поскольку у таких групп любые
 два гомоморфизма
в $\mathbb{Z}$ пропорциональны, найдётся  число $c_\gamma$, для которого
${\rm ind}\chi=c_\gamma{\rm wn}(\chi\circ\gamma)$ при всех $\chi\in X^i$.
Если  в последнем равенстве положить $\chi=\chi_1$, где   $\chi_1\in X^i,\ {\rm ind}\chi_1=1$,
то получим $c_\gamma=k_\gamma^{-1}$,
где $k_\gamma={\rm wn}(\chi_1\circ\gamma)$, что и требовалось доказать.

 ЗАМЕЧАНИЕ 1. Утверждение {\rm 2)} этой теоремы показывает, что в случае линейного
 порядка группа $X^i$ либо порядково
 изоморфна  $\Bbb{Z}$, либо тривиальна.

ЗАМЕЧАНИЕ 2. Из утверждений {\rm 3)} и  {\rm 4)} следует, что определённый выше индекс
 вращения функций из $\Phi(G)$ обладает
 свойствами классического индекса вращения и, в частности, является гомотопическим инвариантом.

ЗАМЕЧАНИЕ 3. Характер $\chi\in T$ называется положительным конечным элементом группы $X$, если для
любого $\xi\in T\setminus\{\bf{1}\}$ найдётся такое натуральное $n$, что
$\xi^n\geq\chi$. Группа $F(X)$, порождённая всеми положительными конечными элементами,
называется группой конечных элементов группы $X$ (эта группа играет важную роль
в \cite{Mur93}). Из теоремы 2 легко следует, что $X^i=F(X)$, если $X^i$ нетривиальна.
В то же время рассмотрение подгрупп $X\subseteq \mathbb{R}$ (наделённых естественным
порядком и дискретной топологией) показывает, что $F(X)$ может быть нетривиальной при
тривиальной $X^i$.

%\begin{corollary}
\textbf{Следствие 1.} \textit{Если группа  $X$  л. у. и  совпадает с $X^i$, то $G$ изоморфна
одномерному тору $\mathbb{T}$}.
%\end{corollary}

В самом деле, так как группа $X$  предполагается нетривиальной,
это сразу следует из замечания 1.

Для случая одномерного тора известно, что все {\it полукоммутаторы}
$[T_\varphi, T_\psi):=T_\varphi T_\psi-T_{\varphi \psi}\ (\varphi,
\psi\in C(\mathbb{T}))$ компактны в $H^p(\mathbb{T})$. Для групп,
отличных от $\mathbb{T}$, это уже не так.

%\begin{corollary}
\textbf{Следствие 2.} \textit{Пусть группа  $X$ л. у., $\chi\in T$. Полукоммутатор $[T_\chi,
T_{\bar{\chi}})$ компактен в $H^p(G)$ тогда и только тогда,
когда $\chi\in X^i$. Следовательно, если операторы $[T_\chi,
T_{\bar{\chi}})$ компактны в $H^p(G)$ при всех $\chi\in T$, то $G$
изоморфна одномерному тору.
}%\end{corollary}

Доказательство. Пусть $I$ --- единичный оператор в $H^p(G)$.
Поскольку $(I-T_\chi T_{\bar{\chi}})\zeta=0$ при
$\zeta\in \chi T$ и  $(I-T_\chi T_{\bar{\chi}})\zeta=\zeta$ при
$\zeta\in T\setminus\chi T$, а $T$ --- линейно независимая
система, порождающая
пространство $H^p(G)$ (лемма 1), оператор $-[T_\chi,
T_{\bar{\chi}})=I-T_\chi T_{\bar{\chi}}$
является ограниченным проектором на подпространство $L \subset H^p(G)$,
порождённое  системой $T\setminus\chi T$. Компактность этого
проектора равносильна конечномерности $L$, т. е. конечности множества
$T\setminus\chi T$, что доказывает первое утверждение. Если теперь операторы $[T_\chi,
T_{\bar{\chi}})$ компактны в $H^p(G)$ при всех $\chi\in T$, то
 $T\subset X^i$, а потому $X=X^i$ и осталось воспользоваться
следствием 1.

%\begin{corollary}
\textbf{Следствие 3.}  \textit{Если группа  $X\not=\mathbb{Z}$ линейно и архимедово
упорядочена, то  $X^i$ тривиальна.
}%\end{corollary}

Действительно, известно, что линейно и архимедово  упорядоченная
группа не содержит собственных нетривиальных выпуклых подгрупп, а равенство
$X=X^i$ влечёт $X=\mathbb{Z}$ в силу утверждения 2) теоремы 2.

\textbf{Предложение 1.} %\begin{proposition}
\textit{Пусть   $\varphi, \psi\in L^\infty(G)$.   Если  $\bar{\varphi}\in
H^\infty(G)$ или $\psi\in H^\infty(G)$, то
}$T_{\varphi\psi}=T_{\varphi}T_\psi$.
%\end{proposition}

Доказательство. 1) Если  $\psi\in H^\infty(G)$, то для любого
$f\in  H^p(G)$
$$
 T_{\varphi}T_\psi  f=T_{\varphi}(\psi f)=P_T(\varphi\psi
 f)=T_{\varphi\psi}f.
$$

2) Пусть теперь  $\bar{\varphi}\in H^\infty(G)$. В силу
доказанного выше при любых $\chi_1, \chi_2\in T$ справедливо
равенство
$$
\langle T_{\bar{\varphi}\psi}\chi_1, \chi_2\rangle=
\langle T_\psi T_{\bar{\varphi}}\chi_1, \chi_2\rangle.
\eqno(3.1)
$$
Как и в классическом случае $G=\mathbb{T}$, легко проверить, что
при $f, g \in H^2(G),\ \varphi\in L^\infty (G)$
$$
\langle T_{\varphi}f,g\rangle = \langle f,T_{\bar{\varphi}}g\rangle.\eqno(3.2)
$$
С учётом последнего равенства левая часть (3.1) приобретает вид
 $\langle \chi_1, T_{\varphi\bar{\psi}}\chi_2\rangle$.

С другой стороны, дважды применяя   равенство (3.2), для правой
части (3.1) получаем выражение $\langle \chi_1,
T_{\varphi}T_{\bar{\psi}}\chi_2\rangle$.
 Следовательно, полагая $h:=T_{\varphi\bar{\psi}}\chi_2-T_{\varphi}T_{\bar{\psi}}\chi_2$,
 выводим, что
$0=\langle \chi_1,h \rangle=\overline{\widehat h(\chi_1)}$, т. е. ${\widehat h}|T=0$.
Но $h\in H^p(G)$, а потому ${\widehat h}|(X\setminus T)=0$. Таким образом,
$T_{\varphi\bar{\psi}}\chi_2=T_{\varphi}T_{\bar{\psi}}\chi_2$  при всех $\chi_2\in T$,
т. е.  $T_{\varphi\bar{\psi}}|Pol_T(G)=T_{\varphi}T_{\bar{\psi}}|Pol_T(G)$, откуда и следует
утверждение предложения.

% \begin{proposition}
\textbf{Предложение 2.}  \textit{Пусть группа  $X$ л. у. Для любой функции $\varphi\in C(G)$  оператор $T_{e^\varphi}$
 обратим в $H^p(G)\ (1<p<\infty)$} .
% \end{proposition}

Доказательство. Воспользуемся методом доказательства теоремы 7.1
из \cite{MurIEOT92}. По теореме Вейерштрасса-Стоуна найдётся такой
тригонометрический полином $q\in Pol(G)$, что
$\|{\bf 1}-e^{\varphi-q}\|_\infty<1/c_p\ (c_p=\|P_T\|_p)$. Пусть
$q_1=P_Tq,\ q_2=q-q_1$. Тогда $q_1, \overline{q_2}\in Pol_T(G)$ (у
нас $X=T\cup T^{-1}$). Лемма 2 показывает теперь, что $e^{\pm
q_1}, e^{\pm\overline{q_2}}\in H^\infty(G)$, а потому в силу
предложения 1 существуют обратные операторы
$(T_{e^{q_i}})^{-1}=T_{e^{-q_i}}\ (i=1,2)$. Повторное применение
предложения 1 даёт $T_{e^\varphi}=T_{e^{q_2}}T_{e^{\varphi
-q}}T_{e^{q_1}}$. Осталось заметить, что оператор $T_{e^{\varphi
-q}}$ обратим, поскольку
$$
\|I-T_{e^{\varphi -q}}\|=\|T_{{\bf 1}-e^{\varphi -q}}\|\leq
c_p\|{\bf 1}-e^{\varphi -q}\|_\infty<1.
$$

Для обобщения теоремы Хартмана-Винтнера
(относительно классического случая см., например, \cite{BS}) нам потребуется следующий
простой факт. Ниже через $w$-$\lim_{\chi\in T}h_\chi$
($s$-$\lim_{\chi\in T}h_\chi$)  обозначается слабый
(соответственно, сильный) предел направленности  $h:T\to
LB(L^p(G))$ ($LB(L^p(G))$ обозначает алгебру линейных ограниченных  операторов в $L^p(G)$).

\textbf{Лемма 4.}  %\begin{lemma}
$w$-$\lim\limits_{\chi\in T}M_\chi=O$.
%\end{lemma}

Доказательство. Пусть $f\in L^p(G),\ g\in L^q(G)$ ($p$ и $q$  --
сопряжённые показатели, $1< p<\infty$). Тогда функция на $X$
$$
\chi\mapsto \langle M_\chi f,g \rangle= \overline{\int\limits_G
\overline\chi \overline f g dm}
$$
стремится к нулю на бесконечности вместе с преобразованием Фурье
функции $\overline fg\in L^1(G)$. Следовательно, для любого
$\varepsilon >0$ найдётся  такое конечное множество
$F=\{\xi_1,\ldots,\xi_n\}\subset X$, что $|\langle M_\chi f,g
\rangle|<\varepsilon$ при $\chi\notin F$. Выберем характеры
$\chi_i,\eta_i\in T$ таким образом, что
$\xi_i^{-1}=\chi_i^{-1}\eta_i\ (i=1,\ldots,n)$, и положим
$\chi_\varepsilon=\chi_1 \ldots \chi_n$. Тогда
$\chi_\varepsilon\xi_i^{-1}\in  T$, т. е. $\chi_\varepsilon\geq \xi_i$ при всех $\xi_i\in F$.
Поэтому из неравенства $\chi>\chi_\varepsilon$ следует, что
$\chi\notin F$, а потому и $|\langle M_\chi f,g
\rangle|<\varepsilon$, что и завершает доказательство леммы.

Теперь обобщение теоремы Хартмана-Винтнера доказывается по
существу так же, как и в классическом случае. Приведём доказательство для
полноты изложения (напомним, что для ограниченного оператора $T$ в банаховом пространстве $Y$
запись  $T\in\Phi_+(Y)$ означает,   что его образ ${\rm Im}T$ замкнут, а ядро
 ${\rm Ker}T$ конечномерно, а запись  $T\in\Phi_-(Y)$ означает конечномерность факторпространства
 $Y/{\rm Im T}$; операторы из $\Phi_-(Y)\cup\Phi_+(Y)$ называются полуфредгольмовыми,
 а операторы из $\Phi_-(Y)\cap\Phi_+(Y)$ --- фредгольмовыми в $Y$).

%\begin{theorem}
\textbf{Теорема 3.} \textit{Если оператор $T_\varphi$  полуфредгольмов в $H^p(G) (1<p<\infty)$, то его символ $\varphi$ обратим в
алгебре $L^\infty (G)$.}
%\end{theorem}

Доказательство. Если $T_\varphi\in\Phi_+(H^p(G))$, то обозначим через $K$
проектор $H^p(G)\to {\rm Ker}T_\varphi$.    Тогда
при некотором $\delta>0$ и всех $f\in H^p(G)$ справедливо неравенство
(см., например, \cite[утверждение 1.12 (g)]{BS})
$$
\|T_\varphi f\|+\|Kf\|\geq\delta\|f\|.
$$
Полагая здесь $f=P_Tg$, имеем
$$
\|P_T M_\varphi P_Tg\|+\|KP_Tg\|\geq\delta\|P_Tg\|.
$$
Далее, $g=P_Tg+Qg$, где  $Q=I-P_T$. Поэтому $\|P_Tg\|\geq\|g\|-\|Qg\|$.
С учётом этого, предыдущее неравенство влечёт
$$
\|P_T M_\varphi P_Tg\|+\|P_TKP_Tg\|\geq\delta (\|g\|-\|Qg\|),
$$
то есть
$$
\|P_T M_\varphi P_Tg\|+\|P_TKP_Tg\|+\delta\|Qg\|\geq\delta \|g\|.
$$
Заменяя тут $g$ на $M_\chi g$, получаем в силу изометричности  $M_\chi$, что
$$
\|M_{\chi^{-1}}P_T M_\varphi P_TM_\chi g\|+\|P_TKP_TM_\chi g\|+\delta\|M_{\chi^{-1}}QM_\chi g\|
\geq\delta \|g\|.\eqno(3.3)
$$

Заметим теперь, что  семейство $E_\chi:=M_{\chi^{-1}}P_TM_\chi\ (\chi\in T)$
равномерно ограничено.
Кроме того,    $E_\chi\xi=\xi (\xi\in X)$ при $\chi\xi\in T$  (т. е.
при $\chi\geq \xi^{-1}$), а потому
$\lim_{\chi\in T}E_\chi q=q$ при $q\in Pol(G)$.
Следовательно,  $s$-$\lim_{\chi\in T}E_\chi=I$ и
$s$-$\lim_{\chi\in T}M_{\chi^{-1}}QM_\chi =O$.

В свою очередь отсюда следует, что
$$
 s\mbox{-}\lim\limits_{\chi\in T}M_{\chi^{-1}}P_T M_\varphi P_TM_\chi=
 s\mbox{-}\lim\limits_{\chi\in T}E_\chi M_\varphi E_\chi=M_\varphi.
 $$
В самом деле,  так как    $E_\chi\xi=\xi$ при $\chi\geq \xi^{-1}$, то   направленность
  $E_\chi M_\varphi E_\chi\xi \to M_\varphi \xi$ по $\chi\in T$.
  По линейности это же верно при
  замене характера  $\xi$ любым тригонометрическим полиномом, и осталось заметить, что
  семейство $E_\chi M_\varphi E_\chi\ (\chi\in T)$   равномерно ограничено.

 Поскольку оператор $P_TKP_T$ компактен, $s$-$\lim_{\chi\in T}P_TKP_TM_\chi =O$
 в силу  леммы 4.  Теперь (3.3) влечёт
 $\|M_\varphi g\|\geq\delta \|g\|$, откуда
 $\|\varphi\|_\infty\geq\delta$.

Пусть теперь $T_\varphi\in\Phi_-(H^p(G))$. Положим ${\stackrel{\circ}{H^p_-}}(G)={\rm Ker} P_T$.
Тогда $L^p(G)=H^p(G)\dot + {\stackrel{\circ}{H^p_-}}(G),\ L^q(G)=H^q(G)\dot +
{\stackrel{\circ}{H^q_-}}(G)$ ($p$ и $q$  --
сопряжённые показатели), а потому
$$
T_\varphi\in\Phi_-(H^p(G))\Longleftrightarrow T_\varphi P_T+Q\in\Phi_-(L^p(G)),
$$
$$
T_{\bar{\varphi}}\in\Phi_+(H^q(G))\Longleftrightarrow T_{\bar{\varphi}} P_T+Q\in\Phi_+(L^q(G)).
$$
Поскольку пространство $L^q(G)$ сопряжено $L^p(G)$, а оператор $T_{\bar{\varphi}} P_T+Q$
сопряжён оператору $T_\varphi P_T+Q$ (см. (3.2)), то $T_{\bar{\varphi}}\in\Phi_+(H^q(G))$, и
всё свелось к случаю, рассмотренному выше.

Сейчас мы в состоянии доказать теорему об индексе.

%\begin{theorem}
\textbf{Теорема 4.}  \textit{Пусть группа  $X$ л. у., $\varphi\in C(G)$. Оператор
$T_\varphi$ в пространстве $H^p(G)\ (1<p<\infty)$ фредгольмов
тогда и только тогда, когда $\varphi\in \Phi(G)$. При этом}
$$
{\rm Ind} T_{\varphi}=-{\rm ind}\varphi.
$$
%\end{theorem}

Доказательство. Необходимость. Если  оператор $T_\varphi$
фредгольмов, то $\varphi\in C(G)^{-1}$ по предыдущей теореме.
Пусть $\varphi=\chi e^g\ (g\in C(G),\ \chi\in X)$ --  разложение
Бора-ван Кампена. По предложению 1 $T_\varphi=T_{e^g}T_\chi$, если
$\chi\in T$, и $T_\varphi=T_\chi T_{e^g}$, если $\chi\in T^{-1}$.
При этом оператор  $T_{e^g}$ обратим  (предложение 2), а потому
фредгольмов нулевого индекса. Следовательно, оператор $T_\chi$
тоже фредгольмов, причём ${\rm Ind} T_{\varphi}={\rm Ind}
T_{\chi}$.

Рассмотрим вопрос о фредгольмовости и индексе оператора $T_\chi$.
 Возможны два случая.

1)  $\chi\in T$. Поскольку ${\rm Ker} T_\chi=\{0\}$,  оператор  $T_\chi$ фредгольмов
тогда и только тогда, когда пространство ${\rm Coker}
T_\chi=H^p(G)/\chi H^p(G)$ конечномерно, и при этом ${\rm Ind}
T_{\chi}=-\dim (H^p(G)/\chi H^p(G))$. Заметим, что соотношения
ортогональности для характеров влекут равенство $\langle\xi,\eta
\rangle=0$ при $\xi\in T\setminus \chi T, \eta\in\chi T$. По
непрерывности $\langle\xi,f \rangle=0$ при $\xi\in T\setminus \chi
T, f\in\chi H^p(G)$. Поэтому характеры из $T\setminus \chi T$
попарно не эквивалентны  ${\rm mod}(\chi H^p(G))$.
 Тогда из сказанного выше следует, что если
$\chi\notin X^i$, т. е. множество $T\setminus \chi T$ бесконечно,
то оператор $T_\chi$ не фредгольмов. Пусть теперь  $\chi\in X^i$,
т. е. $T\setminus \chi T$ конечно. Тогда подпространство ${\rm
span}(T\setminus \chi T)\dot +\chi H^p(G)$ замкнуто в $H^p(G)$
(см., например, \cite[Глава I,  п. 3.3]{Shef}). Поскольку оно
содержит $Pol_T(G)$, справедливо равенство ${\rm span}(T\setminus
\chi T)\dot +\chi H^p(G)=H^p(G)$ (лемма 1), откуда $H^p(G)/\chi H^p(G)={\rm
span}(T\setminus \chi T)$. Следовательно, в рассматриваемом случае
$T_\chi$  фредгольмов и

$$
{\rm Ind} T_{\chi}=-\dim (H^p(G)/\chi H^p(G))=-{\rm ind} \chi=-{\rm
ind}\varphi.
$$

2) Пусть теперь $\chi\in T^{-1}$. Далее мы неоднократно будем
пользоваться равенством $I=T_\chi T_{\bar{\chi}}$, вытекающим из   предложения
1. Если оператор $T_\chi$ фредгольмов, то это равенство
показывает, что оператор $T_{\bar{\chi}}$  тоже фредгольмов. Тогда
$\chi^{-1}\in X^i$ в силу 1), а потому и $\chi\in X^i$. Более
того, из этого же равенства следует, что
$$
{\rm Ind} T_{\chi}=-{\rm Ind} T_{\chi^{-1}}=-{\rm ind} \chi=-{\rm
ind}\varphi.
$$

Обратно, если $\chi\in X^i$, то и $\chi^{-1}\in X^i$, а потому
оператор $T_{\bar{\chi}}$   фредгольмов. Третий раз
воспользовавшись нашим равенством, получаем фредгольмовость
оператора $T_\chi$ и в этом случае.

Итак, в любом случае оператор $T_\chi$ фредгольмов тогда и только
тогда, когда $\chi\in X^i$. Из этого включения, в частности, следует, что
$\varphi\in \Phi(G)$, и доказательство необходимости завершено.
Попутно мы установили и формулу для ${\rm Ind} T_{\varphi}$.

Достаточность. Если $\varphi\in \Phi(G)$, т. е. $\varphi=\chi
e^g$, где  $g\in C(G),\ \chi\in X^i$, то оператор $T_\chi$
фредгольмов по доказанному выше. Тогда оператор $T_{\varphi}$ тоже
фредгольмов  (он равен $T_{e^g}T_\chi$, если $\chi\in T$, и
$T_\chi T_{e^g}$, если $\chi\in T^{-1}$), что и завершает
доказательство.

ЗАМЕЧАНИЕ 4. В работе \cite{Mur93} рассматривалась обобщённая фредгольмовость и
обобщённый индекс (Бройора) тёплицевых операторов в $H^2(G)$ с непрерывным символом
 для случая, когда $X$ счётна и группа её конечных элементов $F(X)$ нетривиальна.
В частности, при этих условиях  в \cite[теорема 4.2]{Mur93} показано, что оператор $T_{\varphi}$
фредгольмов в этом обобщённом смысле тогда и только тогда, когда характер $\chi$
из разложения Бора-ван Кампена $\varphi=\chi e^g$ принадлежит группе $F(X)$.
Сравнение этого результата с теоремой 4 и замечанием 3 показывает, что в случае,
когда группа $X^i$ нетривиальна, обобщённая фредгольмовость оператора $T_{\varphi}$
равносильна  обычной.

Всюду далее  предполагается, что $X$ л. у.

%\begin{corollary}
\textbf{Следствие 4.}  \textit{Пусть  $\varphi\in C(G)$. Следующие утверждения равносильны:}

{\rm 1)} $T_\varphi$ \textit{фредгольмов нулевого индекса};

{\rm 2)}  $\varphi\in \exp (C(G))$;

{\rm 3)} $T_\varphi$ \textit{ обратим}.
%\end{corollary}

%\begin{corollary}
\textbf{Следствие 5.}  \textit{Пусть  $\varphi\in C(G)$, а группа $X\ne
\mathbb{Z}$ линейно и архимедово упорядочена. Если оператор
$T_\varphi$ фредгольмов, то он обратим.}
%\end{corollary}

Доказательство. В силу следствия 3 в этой ситуации имеет место
равенство  $\Phi(G)=\exp(C(G))$.

Справедливость следствия 5 была отмечена в \cite[c. 356]{Mur91}
(там рассматривался  случай  $p=2$).

%\begin{corollary}
\textbf{Следствие 6.} \textit{Если любой оператор  $T_\chi, \ \chi\in T$
фредгольмов, то $G$ изоморфна $\mathbb{T}$}.
%\end{corollary}

Доказательство.  Если    $T_\chi$ фредгольмов, то $\chi\in
\Phi(G)\cap T=X^i\cap T$. Поэтому из условия следует, что
$T\subset X^i$, откуда $X=X^i$ и осталось применить следствие 1.

Рассмотрим примеры вычисления индекса, первый из которых
обобщает  теоремы 3.1 и 3.2 из работы \cite{Mur91}, где рассматривались
пространства $H^2(\Bbb{T}^d)$ и символы специального вида.

{\bf Пример 1}.
  Пусть $G=\Bbb{T}^d,\ d>1$,
$X=\Bbb{Z}^d_{{\rm lex}}$ (нижний индекс указывает, что порядок лексикографический),
 и $\varphi\in C(\Bbb{T}^d)$.
Тогда оператор $T_{\varphi}$ фредгольмов в $H^p(G)\ (1<p<\infty)$, если и только если
его символ имеет вид

$$
\varphi(t_1,\ldots,t_d)= t_d^{n_d}e^{g(t_1,\ldots,t_d)}
$$
 для некоторых  $n_d\in\Bbb{Z},\ g\in C(\Bbb{T}^d)$.

 При этом
$$
{\rm Ind}T_{\varphi}=-n_d.
$$

 В самом деле, по определению лексикографического порядка
 положительный конус в этом случае суть
 $$
   T=\{n\in \mathbb{Z}^d| n_1>0\}\sqcup\{n\in \mathbb{Z}^d| n_1=0, n_2>0\}
   \sqcup
$$
$$
   \cdots\sqcup \{n\in \mathbb{Z}^d| n_1= n_2=\ldots=n_{d-1}=0,
   n_d>0\}\sqcup\{0\}.
 $$
Отсюда следует, что $X^i= \{n\in \mathbb{Z}^d| n_1=
n_2=\ldots=n_{d-1}=0\}$ (точка $(0,\ldots,0,1)$ принадлежит
$X^i$ и осталось воспользоваться замечанием 1 после теоремы 2).
Это значит, если отождествить точку $n\in \Bbb{Z}^d$  с характером
$\chi_n(t)=t^n$ группы  $\Bbb{T}^d$, что
 $$
   \Phi(G)=\{t_d^{n_d}e^{g(t)}| g\in C(\Bbb{T}^d)\},
 $$
 и утверждение следует из  теоремы 4.

{\bf Пример 2}.
  Пусть $G=\Bbb{T}^\infty,$ и
$X=\Bbb{Z}^\infty_{{\rm lex}}$ --- аддитивная группа всех финитных последовательностей целых чисел
с положительным конусом
$$
   T=\{0\}\sqcup\{n\in \mathbb{Z}^\infty| n_1>0\}\sqcup\{n\in \mathbb{Z}^\infty| n_1=0, n_2>0\}
   \sqcup\cdots.
$$

В этом примере группа $X^i$ тривиальна, а потому оператор $T_{\varphi}$
с непрерывным символом фредгольмов в $H^p(G)\ (1<p<\infty)$, если и только если он обратим, т. е.
его символ имеет вид $\varphi=e^g$, где $g\in C(\Bbb{T}^\infty)$.
\footnote[1]{Если  $\varphi\in C(G)$, то оператор $T_\varphi$ по крайней мере односторонне
обратим тогда и только тогда, когда $\varphi\in C(G)^{-1}$. При этом условии он обратим слева
(справа), если и только если характер в разложении Бора-ван Кампена его символа
принадлежит $T^{-1}$ (соответственно, $T$). Это следует из предложений 1, 2 и
теоремы 3. Таким образом, в отличие от классического случая, в описанной ситуации
могут существовать полуфредгольмовы не фредгольмовы тёплицевы операторы с непрерывным
символом. В частности,  так будет в  примере 2. (Это замечание принадлежит рецензенту).}

{\bf Пример 3}.
  Пусть $G=\Bbb{T}^\infty,$
$X=\Bbb{Z}^\infty$, причём линейный порядок в $X$ задаётся положительным конусом
 $T$, состоящим из нуля и тех точек группы $\Bbb{Z}^\infty$, последняя ненулевая координата
 у которых положительна. Поскольку  $(1,0,0,\ldots)\in X^i$, группа $X^i$ состоит из точек
 вида $(n_1,0,0,\ldots)$, где $n_1\in \mathbb{Z}$. Рассуждая как в примере 1, получаем,
 что
 оператор $T_{\varphi}$, где $\varphi\in C(G)$, фредгольмов в $H^p(G)\ (1<p<\infty)$
 тогда и только тогда, когда его символ имеет вид

$$
\varphi(t_1,t_2,\ldots)= t_1^{n_1}e^{g(t_1,t_2,\ldots)}
$$
 для некоторых  $n_1\in\Bbb{Z},\ g\in C(\Bbb{T}^\infty)$. Кроме того,

$$
{\rm Ind}T_{\varphi}=-n_1.
$$

\vspace{5mm}
\centerline{\bf \S 4. Спектры и существенные спектры тёплицевых операторов}
\vspace{5mm}

Применим полученные результаты к исследованию структуры спектров  тёплицевых
операторов. Всюду ниже, если не оговорено противное, группа $X$ считается линейно
упорядоченной, а символ $\varphi$ --  непрерывным.

В случае  символа $\varphi\in C(\mathbb{T})$ известно, что существенный спектр
Фредгольма  $\sigma_e(T_\varphi)=\varphi(\mathbb{T})$. Для других групп ситуация иная.

%  \begin{corollary}
\textbf{Следствие 7.} \textit{Если  $\varphi\in L^\infty(G)$, то
 $\sigma_e(T_\varphi)\supseteq \sigma(\varphi)$, где  $\sigma(\varphi)$ ---
 множество существенных значений
 функции $\varphi$. Причём, если
$G$ не изоморфна $\mathbb{T}$,  для  $\varphi\in
C(G)^{-1}\setminus\Phi(G)$ включение строгое.}
%\end{corollary}

Доказательство. Первое утверждение сразу следует из теоремы 3.
Пусть теперь $\varphi\in C(G)^{-1}\setminus\Phi(G)$ (последнее множество не пусто
по следствию 1, если
$G$ не изоморфна $\mathbb{T}$). Тогда
оператор $T_\varphi$ не фредгольмов (теорема 4), а потому $0
\in\sigma_e(T_\varphi)$.

Далее под {\it дырами} компактного связного множества $\varphi(G)$ подразумеваются  ограниченные компоненты
его дополнения  $\mathbb{C}\setminus\varphi(G)$.

%\begin{theorem}
\textbf{Теорема 5.}   \textit{Пусть  $\varphi\in C(G)$.    Тогда
}
{\rm 1)} \textit{существенный спектр  Фредгольма $\sigma_e(T_\varphi)$ получается из
множества $\varphi(G)$
 присоединением тех его дыр $\Lambda$   (если они существуют), для
которых $\varphi-\lambda\notin \Phi(G)$ при всех (при одном)} $\lambda\in
\Lambda$;

{\rm 2)} \textit{спектр   $\sigma_(T_\varphi)$ получается из
 $\sigma_e(T_\varphi)$ присоединением тех дыр $\Lambda$ множества $\varphi(G)$
 (если они существуют), для
которых $\varphi-\lambda\in \Phi(G)\setminus \exp(C(G))$ при всех (при одном) }$\lambda\in
\Lambda$;

{\rm 3)} \textit{существенный спектр Вейля}
$\sigma_{w}(T_\varphi)=\sigma(T_\varphi)$;

{\rm 4)} \textit{спектры  $\sigma_e(T_\varphi)$ и $\sigma(T_\varphi)$
связны.}
%\end{theorem}

Доказательство. 1) Теорема Бора-ван Кампена
вкупе с единственностью разложения означает, что факторгруппа
$C(G)^{-1}/\exp(C(G))$ изоморфна $X$, поэтому далее мы будем
эти группы отождествлять (значительно более общий факт установлен
в \cite{Gor}).  Обозначим через $\pi$
каноническое отображение
$$
C(G)^{-1}\to C(G)^{-1}/\exp(C(G))=X
$$
 и рассмотрим отображение
$$
\alpha_\varphi: \mathbb{C}\setminus\varphi(G)\to
X, \lambda\mapsto\pi(\varphi-\lambda).
$$
   Оно непрерывно, а значит постоянно на каждой компоненте множества
$\mathbb{C}\setminus\varphi(G)$ в силу дискретности $X$. В соответствии
с теоремой 4
$$
\sigma_e(T_\varphi)=\{\lambda\in\mathbb{C}|  T_{\varphi-\lambda}\notin \Phi(H^p((G))\}=
\{\lambda\in\mathbb{C}|  \varphi-\lambda\notin \Phi(G)\}=
$$
$$
=\varphi(G)\sqcup\{\lambda\in\mathbb{C}
\setminus\varphi(G)|  \alpha_\varphi(\lambda)\notin X^i\}. \eqno(4.1)
$$

Поскольку оператор $T_{\varphi-\lambda}$ обратим при больших $\lambda$,
и $\alpha_\varphi$ постоянно на  компонентах множества
$\mathbb{C}\setminus\varphi(G)$,  неограниченная компонента этого множества
не пересекается с $\sigma_e(T_\varphi)$, а потому формула (4.1) равносильна
первому утверждению теоремы.

2) Воспользуемся равенством
$$
\sigma(T_\varphi)=\varphi(G)\sqcup\{\lambda\in\mathbb{C}\setminus\varphi(G)|
\alpha_\varphi(\lambda)\ne
{\bf 1}\},\eqno(4.2)
$$
которое следует из того, что в силу следствия 4 $\lambda\notin
\sigma(T_\varphi)$ тогда и только тогда, когда $\varphi-\lambda\in
\exp(C(G))$, т. е. когда $\lambda\notin \varphi(G)$ и
$\alpha_\varphi(\lambda)= {\bf 1}$. Осталось заметить, что с учётом формулы (4.1)
 правая часть (4.2) есть
$$
\sigma_e(T_\varphi)\sqcup\{\lambda\in\mathbb{C}\setminus\varphi(G)|
 \alpha_\varphi(\lambda)\in X^i\setminus\{{\bf 1}\}\}.
$$

3) В соответствии с одним утверждением М. Шехтера \cite{Shech} (см.
также \cite[теорема XI.6.12]{Conw}; \cite[ Глава 3]{Er})
$\sigma_{w}(T_\varphi)$ получается из $\sigma_e(T_\varphi)$
присоединением тех точек $\lambda\in \mathbb{C}$, для которых
оператор $T_{\varphi-\lambda}$ полуфредгольмов и ${\rm Ind}T_{\varphi-\lambda}\ne 0$.
С учётом теоремы 4 и утверждения 2) это означает, что $\sigma(T_\varphi)\subseteq\sigma_{w}(T_\varphi)$,
что и требовалось доказать.

4) Это следует из того, что
$\sigma(T_\varphi)$ и $\sigma_e(T_\varphi)$ представимы в виде  объединения связного множества
$\varphi(G)$ и (возможно) некоторых  его дыр.

ЗАМЕЧАНИЕ 5. При $p=2$ связность    $\sigma(T_\varphi)$ была доказана в работе
\cite{MurIrish91}, а связность $\sigma_e(T_\varphi)$ --- в \cite{MurIEOT92}
даже для случая пространств $H^2$, порождённых алгебрами Дирихле.

\textbf{Следствие 8.} %\begin{corollary}
\textit{Спектральные радиусы} $|T_\varphi|_{ess}=|T_\varphi|=\|\varphi\|_\infty$.
%\end{corollary}

Доказательство. Поскольку
 $\sigma(T_\varphi)$ получается из множества $\varphi(G)$ присоединением некоторых
его дыр, то $|T_\varphi|=\sup|\varphi(G)|$. Первое равенство доказывается
аналогично.

\textbf{Следствие 9.} %\begin{corollary}
\textit{Для любого компактного  оператора
 $K$ в пространстве  $H^p(G)\ (1<p<\infty)$ справедливо неравенство
   $\|T_\varphi+K\|\geq (1/c_p)\|T_\varphi\|$, где $c_p=\|P_T\|$.
   В частности, если оператор $T_\varphi$ компактен в  $H^p(G)$,
 то он нулевой.
}%\end{corollary}

Доказательство. Поскольку существенный спектр Вейля инвариантен относительно компактных
возмущений, имеем с учётом теоремы 5,  следствия 8 и теоремы 1
$$
\|T_\varphi+K\|\geq |T_\varphi+K|\geq|T_\varphi|=\|\varphi\|_\infty\geq (1/c_p)\|T_\varphi\|.
$$

\textbf{Следствие 10.}%\begin{corollary}
  \textit{Пусть  $p=2$. Если $A$ есть замкнутая симметричная подалгебра
  алгебры $C(G)$, обладающая тем свойством, что при любых $\varphi, \psi\in A$ полукоммутаторы
  $[T_\varphi, T_\psi)$ компактны в $H^2(G)$, то $\sigma_e(T_\varphi)
=\varphi(G)$ при $\varphi\in A$.}
%\end{corollary}

Доказательство. Отображение $j$ из $A$ в алгебру Калкина ${\cal C}(H^2(G))$,
ставящее в соответствие каждой функции  $\varphi\in A$ класс смежности
$T_\varphi +K(H^2(G))$, является *-гомоморфизмом алгебр, который изометричен
в силу теоремы 1 и  следствия 9. Остаётся заметить, что алгебра $j(A)$, будучи
$C^*$-подалгеброй алгебры Калкина, наполнена в ней.

ЗАМЕЧАНИЕ 6. Примером алгебры $A$, удовлетворяющей условиям следствия 10, может служить
$C^\ast$-подалгебра алгебры $C(G)$, порождённая характером $\chi\in X^i$,
поскольку компактность полукоммутаторов  $[T_{\chi^k},T_{\chi^j}),\ k,j\in\mathbb{Z}$
легко следует из предложения 1 и следствия 2. Следовательно, в этом случае
$\sigma_e(T_\chi)=\chi(G)$. На самом деле, как показывает следствие 11, последнее равенство верно для
любого $p\in (1,\infty)$ (ниже $\mathbb{D}$ обозначает открытый единичный круг в $\mathbb{C}$).

\textbf{Следствие 11.} %\begin{corollary}
  \textit{Пусть  $\chi\in X,\ \chi\ne \bf{1}$.
Тогда}

{\rm 1)} \textit{спектр $\sigma(T_\chi)$ совпадает с замкнутым единичным
кругом $\overline{\mathbb{D}}$};

{\rm 2)} \textit{при $\chi\in X^i$ существенный спектр Фредгольма $\sigma_e(T_\chi)
=\chi(G)=\Bbb{T}$;}

{\rm 3)} \textit{при $\chi\notin X^i$   существенный спектр Фредгольма $\sigma_e
(T_\chi)= \overline{\mathbb{D}}$.
}%\end{corollary}

Доказательство. 1) Равенство  $\chi(G)=\mathbb{T}$
следует из того, что $\chi(G)$ --- нетривиальная связная подгруппа окружности. При
этом для единственной    ограниченной компоненты $\mathbb{D}$
дополнения $\mathbb{C}\setminus \chi(G)$ имеем, используя
обозначения из доказательства предыдущей теоремы,
$\alpha_\chi|\mathbb{D}=\alpha_\chi(0)\ne {\bf 1}$, поскольку
$\chi\notin\exp(C(G))$ в силу единственности разложения Бора-ван
Кампена. Теперь равенство $\sigma(T_\chi)=\overline{\mathbb{D}}$
следует из формулы (4.2).

2) Это вытекает из теоремы 5 и того факта, что в рассматриваемом случае
для единственной дыры $\mathbb{D}$ множества $\chi(G)$ имеем $\mathbb{D}\cap
\sigma_e(T_\chi)=\emptyset$, так как $\chi-0\in \Phi(G)$.

3) Поскольку теперь $\chi-0\notin \Phi(G)$, то $\sigma_e(T_\chi)=\mathbb{T}\cup\mathbb{D}$
по теореме 5.

\

Статья опубликована в \cite{SbMath}.

 \end{document}